\documentclass[12pt]{article}
\usepackage{amsfonts,latexsym,amssymb}
\usepackage{graphicx}
\usepackage[margin=35mm]{geometry}

\newtheorem{LEM}{Lemma}

\newtheorem{THM}{Theorem}
\newtheorem{COR}{Corollary}

\newtheorem{CONJ}{Conjecture}
\newenvironment{Proof}{\textbf{Proof.}}{$\Box$\smallskip}

\begin{document}

\begin{center}
{\Large \textbf{A Note on Bipartite Subgraphs and Triangle-independent Sets}} \medskip \\
Honghai Xu\\
Department of Mathematical Sciences\\
Clemson University\\
Clemson, SC 29634\\
\tt honghax@g.clemson.edu
\end{center}

\begin{abstract}
Let $\alpha_{1} (G)$ denote the maximum size of an edge set that
contains at most one edge from each triangle of $G$.
Let $\tau_{B} (G)$ denote the minimum size of an edge set whose
deletion makes $G$ bipartite.
It was conjectured by Lehel and independently by Puleo that
$\alpha_{1} (G) + \tau_{B} (G) \le n^2/4$
for every $n$-vertex graph $G$.
Puleo showed that $\alpha_{1} (G) + \tau_{B} (G) \le 5n^2/16$
for every $n$-vertex graph $G$.
In this note, we improve the bound by showing that
$\alpha_{1} (G) + \tau_{B} (G) \le 4403n^2/15000$
for every $n$-vertex graph $G$.
\end{abstract}

\noindent \textbf{Keywords:} Bipartite subgraph, Triangle-independent set

\section{Introduction}

Let $G$ be a simple undirected graph.  A \emph{triangle-independent set} in $G$ is
an edge set that contains at most one edge from each triangle of $G$.
We let $\alpha_{1} (G)$ denote the maximum size of a triangle-independent set in $G$.
On the other hand, a \emph{triangle edge cover} in $G$ is an edge set that
contains at least one edge from each triangle of $G$.
We let $\tau_{1} (G)$ denote the minimum size of a triangle edge cover in $G$.

Erd\H{o}s, Gallai, and Tuza made the following conjecture:

\begin{CONJ} (Erd\H{o}s-Gallai-Tuza \cite{erdos1996covering})  \label{cj:EGT}
For every $n$-vertex graph $G$, $\alpha_{1} (G) + \tau_{1} (G) \le n^2/4$.
\end{CONJ}

Note that the equality holds for the graphs $K_n$ and $K_{n/2,n/2}$, where $n$ is even.
Indeed, $\alpha_{1} (K_n) = n/2$ and $\tau_{1} (K_n) = {n \choose 2} - n^2/4$
(by Mantel's theorem~\cite{mantel1907problem}), while
$\alpha_{1} (K_{n/2,n/2}) = n^2/4$ and $\tau_{1} (K_{n/2,n/2}) = 0$.
In both cases, $\alpha_{1} (G) + \tau_{1} (G) = n^2/4$.
More generally, let $G_1 \vee \ldots \vee G_t$ denote the graph obtained from
the disjoint union $G_1 + \ldots + G_t$ by adding all edges between vertices from
different $G_i$. Puleo (see ~\cite{puleo2015conjecture, puleo2015extremal}) showed
that the equality holds for any graph of the form
$K_{r_1,r_1} \vee \ldots \vee K_{r_t,r_t}$.

Conjecture~\ref{cj:EGT} was originally stated only
for \emph{triangular graphs},
which are graphs where every edge lies in a triangle
(see~\cite{erdos1996covering,erdos1990some}).  However, later it was
stated for general graphs (see~\cite{erdos1999selection,tuza2001unsolved}).
It was proved by Puleo~\cite{puleo2015extremal} that these two forms of the
conjecture are equivalent.

A related parameter, denoted by $\tau_{B} (G)$, is the minimum size of
an edge set in $G$ whose deletion makes $G$ bipartite.
Clearly $\tau_{B} (G) \ge \tau_{1} (G)$.  Erd\H{o}s~\cite{erdos1984some} asked
which graphs satisfy $\tau_{B} (G) = \tau_{1} (G)$.  Bondy, Shen, Thomass{\'e}, and
Thomassen~\cite{bondy2006density} proved that $\tau_{B} (G) = \tau_{1} (G)$ when
$\delta (G) \ge 0.85~|V(G)|$, and later Balogh, Keevash, and
Sudakov~\cite{balogh2006minimal} proved that $\tau_{B} (G) = \tau_{1} (G)$ when
$\delta (G) \ge 0.79~|V(G)|$.

The following conjecture, which is stronger than
Conjecture~\ref{cj:EGT}, was proposed by Lehel (see~\cite{erdos1990some}) and
independently by Puleo~\cite{puleo2015conjecture}.

\begin{CONJ} (\cite{puleo2015conjecture}) \label{cj:Puleo}
For every $n$-vertex graph $G$, $\alpha_{1} (G) + \tau_{B} (G) \le n^2/4$.
\end{CONJ}

Puleo~\cite{puleo2015conjecture, puleo2015extremal} obtained many interesting
results towards Conjectures~\ref{cj:EGT} and~\ref{cj:Puleo}.
Conjecture~\ref{cj:Puleo} was verified for triangle-free graphs and for
graphs that have no induced subgraph isomorphic to $K_4^{-}$ (the graph
obtained from $K_4$ by deleting an edge)~\cite{puleo2015extremal}.
For general graphs, Puleo~\cite{puleo2015conjecture} showed the following upper bound:

\begin{THM} (\cite{puleo2015conjecture}) \label{t:PuleoUB}
For every $n$-vertex graph $G$, $\alpha_{1} (G) + \tau_{B} (G) \le 5n^2/16$.
\end{THM}

The main purpose of this note is to provide an improved bound towards
Conjecture~\ref{cj:Puleo}.  We prove that
$\alpha_{1} (G) + \tau_{B} (G) \le 4403n^2/15000$ for every $n$-vertex graph $G$.
We use ideas from~\cite{puleo2015extremal}, \cite{puleo2015conjecture},
~\cite{sudakov2007note}, and~\cite{furedi2015proof}.

We shall use the following notation and terminology.
For shorthand, we let $f_B(G) = \alpha_{1} (G) + \tau_{B} (G)$.
We let $n(G)$, $e(G)$, and $t(G)$ denote the number of vertices, edges, and
triangles in $G$, respectively.  When there is no confusion involved,
we simply write $n$, $e$, and $t$.
We let $d(v)$ denote the degree of a vertex $v$, and $\omega(G)$ denote
the clique number of $G$.  When $S \subseteq V(G)$, we write
$G[S]$ for the subgraph of $G$ induced by $S$, $\overline{S}$
for the set $V(G)-S$, and $[S,\overline{S}]$ for the set of all edges
with one endpoint in $S$ and the other endpoint in $\overline{S}$.
We use the term \emph{minimal counterexample} to
refer to a vertex-minimal counterexample, that is, a graph $G$ such that
the property in question holds for every proper induced subgraph of $G$
but does not hold for $G$.

The rest of the paper is organized as follows.  In the next section,
we investigate the structure of a minimal counterexample to
$f_B(G) \le cn(G)^2$ where $c > 1/4$.
We show that the clique number of such a counterexample is bounded by
a function of $c$.  Thus, to prove that $f_B(G) \le cn(G)^2$, we only need
to prove it for graphs with small clique number.
Then in Section~\ref{s:firstImprovement} we present a quick proof of
$f_B(G) \le 3n(G)^2/10$, which improves Theorem~\ref{t:PuleoUB}.
In Section~\ref{s:K6freeBipartite} we give some estimates of
$\tau_{B} (G)$ for $K_6$-free graphs.  In particular, we show that
every $n$-vertex $K_6$-free graph can be made bipartite by deleting
at most $17n^2/100$ edges.  In Section~\ref{s:Main} we prove our main result.

\section{$f_B(G)$ and clique number} \label{s:clique}

We need the following lemma from~\cite{puleo2015extremal}.

\begin{LEM} (\cite{puleo2015extremal}) \label{l:PuleoDenseCut}
Let $G$ be a graph, and let $A$ be a triangle-independent set of edges in $G$.
If $S$ is a nonempty proper subset of $V(G)$, then
\[
f_B(G) \le f_B(G[S]) + f_B(G[\overline{S}])
+ \frac{1}{2} \left\vert [S,\overline{S}] \right\vert
+ \left\vert [S,\overline{S}] \cap A \right\vert.
\]
\end{LEM}

In~\cite{puleo2015extremal}, Puleo used Lemma~\ref{l:PuleoDenseCut} to
prove some conclusions on the structure of a minimal counterexample $G$
to Conjecture~\ref{cj:Puleo}.  By slightly extending his argument, we
show the following:

\begin{LEM} \label{l:CliqueLemma}
For any constant $c > 1/4$, if $G$ is a minimal counterexample to
$ f_B(G) \le cn(G)^2$, then $\omega(G) < 1/(4c-1)$.
\end{LEM}
\begin{Proof} Let $G$ be a minimal counterexample to $ f_B(G) \le cn(G)^2$.
We may assume $n(G)\ge 5$, since it is easy to verify that
$f_B(G) \le n(G)^2/4 \le cn(G)^2$ when $n(G) \le 4$.
Let $K$ be the largest clique in $G$, and let $k=|K|=\omega(G)$.
Since $f_B(G) \le n(G)^2/4 \le cn(G)^2$ when $G$ is complete, we may
assume $1\le k \le n(G)-1$.

For simplicity, write $n$ for $n(G)$.  Let $A$ be any triangle-independent
set in $G$, and for every $v\in V(G)$, let $N_A(v)=\{\;w\in V(G): vw\in A\;\}$.
Since $A$ is triangle-independent, $\left\vert N_A(v) \cap K \right\vert \le 1$
for each $v \in \overline{K}$.  It follows that
$\left\vert [K,\overline{K}] \cap A \right\vert \le n-k$.
By Lemma~\ref{l:PuleoDenseCut} and the minimality of $G$, we have
\begin{eqnarray*}
cn^2 < f_B(G) &\le& f_B(G[K]) + f_B(G[\overline{K}])
+ \frac{1}{2} \left\vert [K,\overline{K}] \right\vert
+ \left\vert [K,\overline{K}] \cap A \right\vert\\
&\le& \frac{k^2}{4} + c{(n-k)}^2 + \frac{1}{2} \left\vert [K,\overline{K}] \right\vert + n-k.
\end{eqnarray*}
Thus, $\left\vert [K,\overline{K}] \right\vert > -(2c+\frac{1}{2})k^2+4cnk+2k-2n$.
However, since $K$ is the largest clique of $G$, $\left\vert [K,\overline{K}] \right\vert \le (n-k)(k-1)$.
Hence, we have
\begin{eqnarray*}
(n-k)(k-1) > -\left(2c+\frac{1}{2}\right)k^2+4cnk+2k-2n.
\end{eqnarray*}
The above inequality simplifies to $(\frac{1}{2}-2c)k^2 + k < \left(1-(4c-1)k\right)n$.
Assume to the contrary that $k \ge 1/(4c-1)$.  Then
$\left(1-(4c-1)k\right)n \le \left(1-(4c-1)k\right)k$.
It follows that $(\frac{1}{2}-2c)k^2 + k < \left(1-(4c-1)k\right)k$.
That is, $c < 1/4$, a contradiction.
\end{Proof}

\section{A first improvement} \label{s:firstImprovement}

In this section we present a quick proof of $f_B(G) \le 3n(G)^2/10$.
We first show that the conclusion holds for $K_5$-free graphs, and then
use Lemma~\ref{l:CliqueLemma} to prove that it holds for all graphs.

For a graph $G$, let $b(G)$ denote the largest size of a vertex set $B$
such that $B$ induces a bipartite subgraph of $G$.
Puleo~\cite{puleo2015conjecture} proved
the following bound for $\alpha_{1} (G)$:

\begin{LEM} (\cite{puleo2015conjecture}) \label{l:PuleoTriIndependentSet1}
For every $n$-vertex graph $G$, $\alpha_{1} (G) \le nb(G)/4$.
\end{LEM}

Now we consider $\tau_{B} (G)$.  A well-known result by
Erd\H{o}s~\cite{erdos1965some} says that $\tau_{B} (G) \le e/2$ for
every graph $G$ with $e$ edges.  Puleo~\cite{puleo2015conjecture} proved
the following bound for $\tau_{B} (G)$:

\begin{LEM} (\cite{puleo2015conjecture}) \label{l:PuleoMakeBipartite1}
For every $n$-vertex graph $G$, $\tau_{B} (G) \le \left(n^2-b(G)^2\right)/4$.
\end{LEM}

When $G$ is a $K_5$-free graph, $\tau_{B} (G)$ can be bounded as follows:

\begin{LEM} \label{l:K5freeLemma1}
For every $n$-vertex $K_5$-free graph $G$,
\[
\tau_{B} (G) \le \frac{b(G) \left(n-b(G)\right)}{2} + \frac{3{\left(n-b(G)\right)}^2}{16}.
\]
\end{LEM}
\begin{Proof}
Let $B$ denote the vertex set of a largest bipartite induced subgraph of $G$.
Since $G[\overline{B}]$ is $K_5$-free, by
Tur\'{a}n's theorem~\cite{turan1941extremalaufgabe} it has at most
$3{(n-b(G))}^2/8$ edges.
Therefore $G[\overline{B}]$ can be made bipartite by deleting at most
$3{(n-b(G))}^2/16$ edges.
The conclusion follows by considering the two different ways to join
the partite sets of a largest bipartite subgraph in $G[\overline{B}]$ with
the partite sets of $G[B]$.
\end{Proof}

Now we can give the following bound for $f_B(G)$ when $G$ is $K_5$-free.

\begin{THM} \label{t:K5free1}
For every $n$-vertex $K_5$-free graph $G$, $f_B(G) \le 3n^2/10$.
\end{THM}
\begin{Proof}
By Lemma~\ref{l:PuleoTriIndependentSet1} and
Lemma~\ref{l:K5freeLemma1}, we have
\begin{eqnarray*}
f_B(G) = \alpha_{1} (G) + \tau_{B} (G)
&\le& \frac{nb(G)}{4} + \frac{b(G) (n-b(G))}{2} + \frac{3{\left(n-b(G)\right)}^2}{16} \\
&=& \frac{1}{16} \left( -5b(G)^2 + 6nb(G) +3n^2 \right) \\
&=& g(b(G)),
\end{eqnarray*}
where $g(x) = (-5x^2+6nx+3n^2)/16$.  Since $g(x)$ achieves its
maximum at $x = 3n/5$, we have $f_B(G) \le g(3n/5) = 3n^2/10$.
\end{Proof}

By using Theorem~\ref{t:K5free1} and Lemma~\ref{l:CliqueLemma},
we show $f_B(G) \le 3n^2/10$, which improves Theorem~\ref{t:PuleoUB}.

\begin{THM} \label{t:Main1}
For every $n$-vertex graph $G$, $f_B(G) \le 3n^2/10$.
\end{THM}
\begin{Proof} It is easy to verify the conclusion for small $n$.
Assume to the contrary that $G$ is a minimal counterexample.
By Theorem~\ref{t:K5free1}, $\omega (G) \ge 5$.
However, by Lemma~\ref{l:CliqueLemma} we have
$\omega (G) < 1/(4\times \frac{3}{10} - 1) = 5$, a contradiction.
\end{Proof}

\section{$\tau_{B} (G)$ for $K_6$-free graphs} \label{s:K6freeBipartite}

To improve our bound for $f_{B} (G)$, we consider $\tau_{B} (G)$ for
$K_6$-free graphs.  Similar questions have been investigated by
various researchers.  Erd\H{o}s~\cite{erdos1975problems} conjectured that
every $n$-vertex triangle-free graph can be made bipartite by
deleting at most $n^2/25$ edges.  Erd\H{o}s, Faudree, Pach and
Spencer~\cite{erdos1988make} proved that it is enough to delete
$(1/18-\epsilon)n^2$ edges to make a $n$-vertex triangle-free graph bipartite.
Erd\H{o}s (see e.g., ~\cite{erdos1988make}) also conjectured that it is enough
to delete at most $(1+o(1))n^2/9$ edges to make any $n$-vertex $K_4$-free
graph bipartite. This was confirmed by Sudakov~\cite{sudakov2007note}
in the following strong form:

\begin{THM} (\cite{sudakov2007note}) \label{t:Sudakov}
Every $n$-vertex $K_4$-free graph can be made bipartite
by deleting at most $ n^2/9$ edges.  Moreover, equality holds if and
only if $G$ is a complete $3$-partite graph with parts of size $n/3$.
\end{THM}

Furthermore, Sudakov~\cite{sudakov2007note} made the following conjecture
on $\tau_{B} (G)$ for $K_r$-free graphs where $r\ge 5$:

\begin{CONJ} (\cite{sudakov2007note}) \label{cj:SudakovBipartiteConj}
Let $G$ be a $n$-vertex $K_r$-free graph where $r\ge 5$.  Then
\[
\tau_{B} (G) \le
\left\{
\begin{array}{rl}
    \frac{r-3}{4(r-1)}n^2, & \mbox{if $r$ is odd} \\[10pt]
    \frac{{(r-2)}^2}{4{(r-1)}^2}n^2, & \mbox{if $r$ is even}
\end{array}
\right.
\]
\end{CONJ}

This conjecture seems to be very difficult.  The original paper
of Sudakov~\cite{sudakov2007note} pointed out that
some of the ideas there can be used to make a progress on
the conjecture for even $r$.

Our focus in this section is to give some estimates on
$\tau_{B} (G)$ for $K_6$-free graphs.  We first consider bounds on
$\tau_{B} (G)$ for $K_5$-free graphs, and then use the bounds that
we obtain to prove bounds on $\tau_{B} (G)$ for $K_6$-free graphs.
The key ideas that we use come from~\cite{sudakov2007note} and~\cite{furedi2015proof}.
We start with the following well-known fact.

\begin{LEM} (see, e.g., Lemma 2.1 of~\cite{alon1996bipartite}) \label{l:AlonSudakov4Partite}
Let $G$ be a (at most) $2m$-partite graph with $e$ edges.
Then $\tau_{B} (G) \le (m-1)e/(2m-1)$.
\end{LEM}

We also need the following theorem from~\cite{furedi2015proof}, which
is a sharpening of Tur\'{a}n's theorem.
It helps us to deal with the case that the graph is dense:

\begin{THM} (\cite{furedi2015proof}) \label{t:furedi}
Every $n$-vertex $K_{p+1}$-free graph $G$ with $e(T_{n,p})-k$ edges
contains a (at most) $p$-partite subgraph with at least $e(G)-k$ edges, where
$T_{n,p}$ is the complete $p$-partite graph of order $n$ having
the maximum number of edges.
\end{THM}

\begin{COR} \label{c:furediBipartite}
Let $G$ be a graph on $n$ vertices with $e$ edges. \\
(a) If $G$ is $K_5$-free, then $\tau_{B} (G) \le n^2/4 - e/3$; \\
(b) If $G$ is $K_6$-free, then $\tau_{B} (G) \le 6n^2/25 - e/5$.
\end{COR}
\begin{Proof}
Suppose $G$ is $K_5$-free.  Let $H$ be a $4$-partite subgraph
of $G$ having the maximum number of edges.
By Theorem~\ref{t:furedi}, $e(H) \ge 2e - 3n^2/8$.
By Lemma~\ref{l:AlonSudakov4Partite},
$H$ can be made bipartite by deleting at most $e(H)/3$ edges.
Thus,
\[
\tau_{B} (G) \le e-e(H)+\frac{e(H)}{3} = e-\frac{2e(H)}{3}
\le e-\frac{2}{3}\left(2e - \frac{3n^2}{8}\right)
= \frac{n^2}{4} - \frac{e}{3}.
\]
This proves $(a)$.

Suppose $G$ is $K_6$-free.  Let $H$ be a $5$-partite subgraph
of $G$ having the maximum number of edges.
By Theorem~\ref{t:furedi}, $e(H) \ge 2e - 2n^2/5$.
By Lemma~\ref{l:AlonSudakov4Partite},
$H$ can be made bipartite by deleting at most $2e(H)/5$ edges.
Thus,
\[
\tau_{B} (G) \le e-e(H)+\frac{2e(H)}{5} = e-\frac{3e(H)}{5}
\le e-\frac{3}{5}\left(2e - \frac{2n^2}{5}\right)
= \frac{6n^2}{25} - \frac{e}{5}.
\]
This proves $(b)$.
\end{Proof}

\begin{LEM} (see Lemma 2.3 of \cite{sudakov2007note}) \label{l:SudakovBipartite}
Let $G$ be a graph on $n$ vertices with $e$ edges and $t$ triangles.
Then $\tau_{B} (G) \le e + \left(6t-\sum\limits_{v}d^2(v)\right)/n$.
\end{LEM}

Our next step is to apply some of the ideas and techniques from~\cite{sudakov2007note}
to prove a bound on $\tau_{B} (G)$ for $K_5$-free graphs.

\begin{LEM} \label{l:K5freeBipartite1}
Let $G$ be a $K_5$-free graph on $n$ vertices with $e$ edges
and $t$ triangles.  Then
$\tau_{B} (G) \le e/2 + \left(2\sum\limits_{v}d^2(v) - 27t \right)/\left(18n\right)$.
\end{LEM}
\begin{Proof}
Let $v$ be a vertex of $G$ and let $e_v$ denote the number of
edges spanned by the neighborhood $N(v)$.  The induced subgraph
$G[N(v)]$ is $K_4$-free, since $G$ is $K_5$-free.  By Theorem~\ref{t:Sudakov},
$G[N(v)]$ can be made bipartite by deleting
at most $d^2(v)/9$ edges.  Let $A$ and $B$ be the resulting partite sets
of $G[N(v)]$.  We obtain a bipartite
subgraph of $G$ by placing the vertices in $G-N(v)$ into the partite
sets $A$ and $B$ randomly and independently with probability $1/2$, and
deleting all edges within the partite sets.  For each edge
in $G-G[N(v)]$, it is deleted with probability $1/2$.
By linearity of expectation, $\tau_{B} (G) \le (e - e_v)/2 + d^2(v)/9$.
By averaging over all vertices $v$,
we have
\begin{eqnarray*}
\tau_{B} (G) \le \frac{e}{2} + \frac{1}{9n}\sum\limits_{v}d^2(v) - \frac{1}{2n}\sum\limits_{v}e_v
             = \frac{e}{2} + \frac{1}{18n}\left(2\sum\limits_{v}d^2(v) - 27t\right),
\end{eqnarray*}
where we have used the fact that $\sum\limits_{v}e_v = 3t$.
\end{Proof}

Now we can bound $\tau_{B} (G)$ for $K_5$-free graphs
in terms of $n(G)$ only.

\begin{THM} \label{t:K5freeBipartite2}
For every $n$-vertex $K_5$-free graph $G$, $\tau_{B} (G) \le 29n^2/200$.
\end{THM}
\begin{Proof}
By Lemma~\ref{l:SudakovBipartite} and Lemma~\ref{l:K5freeBipartite1},
we have $\tau_{B} (G) \le e + \left(6t-\sum\limits_{v}d^2(v)\right)/n$ and
$\tau_{B} (G) \le e/2 + (2\sum\limits_{v}d^2(v) - 27t)/(18n)$.
Multiplying the first inequality by $1/5$ and the second
inequality by $4/5$, and adding them together, we have
\begin{eqnarray*}
\tau_{B} (G) &\le& \frac{3e}{5} - \frac{1}{9n}\sum\limits_{v}d^2(v) \\
             &\le& \frac{3e}{5} - \frac{1}{9n^2} {\left(\sum\limits_{v}d(v)\right)}^2 \\
             &\le& \frac{3e}{5} - \frac{4e^2}{9n^2},
\end{eqnarray*}
where we have used the Cauchy-Schwartz inequality and
the fact that $\sum\limits_{v}d_v = 2e$.

First consider the case when $e \le 63n^2/200$.  Note that the
function $g(x)=3x/5 - 4x^2/9$ is increasing in the interval
$x \le 63/200$.  So we have
\begin{eqnarray*}
\tau_{B} (G) \le g(e/n^2)n^2 \le g(63/200)n^2 = \frac{1449n^2}{10000} < \frac{29n^2}{200}.
\end{eqnarray*}

Next consider the case when $e > 63n^2/200$.
By Corollary~\ref{c:furediBipartite} $(a)$ we have
\begin{eqnarray*}
\tau_{B} (G) \le \frac{n^2}{4} - \frac{e}{3} <
\frac{n^2}{4} - \frac{21n^2}{200} = \frac{29n^2}{200}.
\end{eqnarray*}
\end{Proof}

\noindent \textbf{Remark.} Since a $K_5$-free graph has at most $3n^2/8$ edges,
it is enough to delete at most $3n^2/16$ edges to make it bipartite.
Although the bound in Theorem~\ref{t:K5freeBipartite2} is better than that,
it probably can be improved substantially.  Indeed
Conjecture~\ref{cj:SudakovBipartiteConj} says that it suffices to delete
$n^2/8$ edges to make a $K_5$-free graph bipartite.
It seems that some new ideas or tools are needed to improve the estimate above.

Next we use the bounds we obtained to prove bounds on
$\tau_{B} (G)$ for $K_6$-free graphs.  The approach is
nearly identical to that used for $K_5$-free graphs.

\begin{LEM} \label{l:K6freeBipartite1}
Let $G$ be a $K_6$-free graph on $n$ vertices with $e$ edges and $t$ triangles.
Then $\tau_{B} (G) \le e/2 + (29\sum\limits_{v}d^2(v) - 300t)/(200n)$.
\end{LEM}
\begin{Proof}
Let $v$ be a vertex of $G$ and let $e_v$ denote the number of edges spanned
by the neighborhood $N(v)$.  The induced subgraph $G[N(v)]$ is $K_5$-free,
since $G$ is $K_6$-free.  By Theorem~\ref{t:K5freeBipartite2}, $G[N(v)]$ can be made
bipartite by deleting at most $29d^2(v)/200$ edges.  Thus $\tau_{B} (G) \le (e - e_v)/2 + 29d^2(v)/200$.
By averaging over all vertices $v$, we have
\begin{eqnarray*}
\tau_{B} (G) \le
\frac{e}{2} + \frac{29}{200n}\sum\limits_{v}d^2(v) - \frac{1}{2n}\sum\limits_{v}e_v =
\frac{e}{2} + \frac{1}{200n}\left(29\sum\limits_{v}d^2(v) - 300t\right),
\end{eqnarray*}
where we have used the fact $\sum\limits_{v}e_v = 3t$.
\end{Proof}

\begin{THM} \label{t:K6freeBipartite2}
For every $n$-vertex $K_6$-free graph $G$, $\tau_{B} (G) \le 17n^2/100$.
\end{THM}
\begin{Proof}
By Lemma~\ref{l:SudakovBipartite} and Lemma~\ref{l:K6freeBipartite1},
we have $\tau_{B} (G) \le e + \left(6t-\sum\limits_{v}d^2(v)\right)/n$ and
$\tau_{B} (G) \le e/2 + (29\sum\limits_{v}d^2(v) - 300t)/(200n)$.
Multiplying the first inequality by $1/5$ and the second
inequality by $4/5$, and adding them together, we have
\begin{eqnarray*}
\tau_{B} (G) &\le& \frac{3e}{5} - \frac{21}{250n}\sum\limits_{v}d^2(v) \\
             &\le& \frac{3e}{5} - \frac{21}{250n^2} {\left(\sum\limits_{v}d(v)\right)}^2 \\
             &\le& \frac{3e}{5} - \frac{42e^2}{125n^2},
\end{eqnarray*}
where we have used the Cauchy-Schwartz inequality and
the fact that $\sum\limits_{v}d_v = 2e$.

First consider the case when $e \le 35n^2/100$.  Note that the
function $g(x)=3x/5 - 42x^2/125$ is increasing in the interval
$x \le 35/100$.  So we have
\begin{eqnarray*}
\tau_{B} (G) \le g(e/n^2)n^2 \le g(35/100)n^2 = \frac{4221n^2}{25000} < \frac{17n^2}{100}.
\end{eqnarray*}

Next consider the case when $e > 35n^2/100$.
By Corollary~\ref{c:furediBipartite} $(b)$ we have
\begin{eqnarray*}
\tau_{B} (G) \le \frac{6n^2}{25} - \frac{e}{5} <
\frac{6n^2}{25} - \frac{7n^2}{100} = \frac{17n^2}{100}.
\end{eqnarray*}
\end{Proof}

\noindent \textbf{Remark.} The bound above is also probably not tight.
Conjecture~\ref{cj:SudakovBipartiteConj} says that it enough to delete
at most $16n^2/100$ edges to make a $K_6$-free graph bipartite.
Nevertheless, it still suffices for our purpose.

To prove our main result, we also need bounds on $\tau_{B} (G)$ for
$K_6$-free graphs in terms of $n$, $e$, and $b(G)$.
Let $B$ be the vertex set of a largest bipartite induced subgraph of $G$.
By a similar argument to that used in Lemma~\ref{l:K5freeLemma1}, we
have that $\tau_{B}(G) \le b(G)(n-b(G))/2 + 17{(n-b(G))}^2/100$.
However this is a very rough estimate.  Indeed, if
$\left\vert [B,\overline{B}] \right\vert = b(G)\left(n-b(G)\right)$, then
since $G$ is $K_6$-free, $G[\overline{B}]$ cannot have many edges and so
it could be made bipartite by deleting less than $17{(n-b(G))}^2/100$
edges.  To refine our argument, we need the following lemma.

\begin{LEM} \label{l:K6freeBipartite3}
Let $G$ be a $K_6$-free graph on $n$ vertices, and let
$S$ be a vertex set that induces a $K_5$-free subgraph of $G$.
If $|S| \ge 49n/50$, then $\tau_{B} (G) \le 3n^2/20$.
\end{LEM}
\begin{Proof}
First note that
$\tau_{B} (G) \le \tau_{B}(G[S]) + \tau_{B}(G[\overline{S}]) + \frac{1}{2}\left\vert [S,\overline{S}] \right\vert$,
which follows by considering the two different ways to join the partite sets
of a largest bipartite subgraph in $G[S]$ with those of one in $G[\overline{S}]$.
Let $s = |S|$.  Since $G[S]$ is $K_5$-free, by Theorem~\ref{t:K5freeBipartite2}
we have $\tau_{B} (G[S]) \le 29s^2/200$.  Since $G[\overline{S}]$ is
$K_6$-free, by Theorem~\ref{t:K6freeBipartite2}
we have $\tau_{B} (G[\overline{S}]) \le 17{(n-s)}^2/100$.
Thus,
\begin{eqnarray*}
\tau_{B} (G) &\le& \frac{29s^2}{200} + \frac{17{(n-s)}^2}{100} + \frac{s(n-s)}{2}\\
              &=& \frac{1}{200}\left(-37s^2 + 32ns + 34n^2\right).
\end{eqnarray*}
The function $g(s) = (-37s^2 + 32ns + 34n^2)/200$ is decreasing
in the interval $s \ge 49n/50$.  Thus, if $s \ge 49n/50$, then
$\tau_{B} (G) \le g(49n/50) = 74563 n^2/500000 < 3n^2/20$.
\end{Proof}

We finish this section by the following corollary.

\begin{COR} \label{c:K6freeBipartite4}
Let $G$ be a $K_6$-free graph on $n$ vertices with $e$ edges.  Then \\
(a)
\begin{eqnarray*}
\tau_{B} (G) \le \max \left(\frac{-7b(G)^2+4nb(G)+3n^2}{20},
\frac{-32b(G)^2+15nb(G)+17n^2}{100} \right).
\end{eqnarray*}
(b)
\begin{eqnarray*}
\tau_{B} (G) \le \frac{e}{3} + \frac{17{(n-b(G))}^2}{150}.
\end{eqnarray*}
\end{COR}
\begin{Proof}
Let $B$ be the vertex set of a largest bipartite induced subgraph of $G$.

We first prove $(a)$.  Note that by considering the two different ways to
join the partite sets of a largest bipartite subgraph in $G[\overline{B}]$ with
the partite sets of $G[B]$, we have
$\tau_{B}(G) \le \tau_{B}(G[\overline{B}]) + \frac{1}{2} \left\vert [B,\overline{B}] \right\vert$.
There are two possible cases:

If there is a vertex $v \in B$ that has at least
$49(n-b(G))/50$ neighbors in $G[\overline{B}]$, then since $G$ is
$K_6$-free, those neighbors of $v$ in $G[\overline{B}]$ must
induce a $K_5$-free subgraph of $G[\overline{B}]$.
By Lemma~\ref{l:K6freeBipartite3}, $G[\overline{B}]$ can be
made bipartite by deleting at most $3{(n-b(G))}^2/20$ edges.
Thus, we have
\begin{eqnarray*}
\tau_{B}(G) \le \frac{3{(n-b(G))}^2}{20} + \frac{b(G)(n-b(G))}{2}
       = \frac{-7b(G)^2+4nb(G)+3n^2}{20}.
\end{eqnarray*}

If every vertex $v \in B$ has at most $49(n-b(G))/50$ neighbors
in $G[\overline{B}]$, then
$\left\vert [B,\overline{B}] \right\vert < 49b(G)(n-b(G))/50$.
Since $G[\overline{B}]$ is $K_6$-free, by
Theorem~\ref{t:K6freeBipartite2} we have
$\tau_{B}(G[\overline{B}]) \le 17{\left(n-b(G)\right)}^2/100$.
It follows that
\begin{eqnarray*}
\tau_{B}(G) \le \frac{17{(n-b(G))}^2}{100} + \frac{49b(G)(n-b(G))}{100}
       =  \frac{-32b(G)^2+15nb(G)+17n^2}{100}.
\end{eqnarray*}
This proves $(a)$.

We next prove $(b)$.  Note that we can make $G$ $4$-partite
by deleting $\tau_{B} (G[\overline{B}])$ edges in $G[\overline{B}]$.
By Lemma~\ref{l:AlonSudakov4Partite}, we can make the resulting $4$-partite graph bipartite
by deleting at most $\left(e-\tau_{B} (G[\overline{B}])\right)/3$ edges.
It follows that
\begin{eqnarray*}
\tau_{B} (G) &\le& \tau_{B} (G[\overline{B}]) + \frac{e-\tau_{B}(G[\overline{B}])}{3}\\
             &\le& \frac{e}{3} + \frac{2\tau_{B} (G[\overline{B}])}{3}\\
             &\le& \frac{e}{3} + \frac{17{(n-b(G))}^2}{150}.
\end{eqnarray*}
This proves $(b)$.
\end{Proof}

\section{Main Result} \label{s:Main}

In this section we prove $f_B(G) \le 4403n(G)^2/15000$.
We first show that the conclusion holds for $K_6$-free graphs, and
then use Lemma~\ref{l:CliqueLemma} to prove that it holds for all graphs.

We need the following lemma from~\cite{puleo2015conjecture}.

\begin{LEM} (\cite{puleo2015conjecture}) \label{l:PuleoTriIndependentSet2}
For every graph $G$ on $n$ vertices with $e$ edges,
$\alpha_{1} (G) \le n^2/2 - e$.
\end{LEM}

\begin{THM} \label{t:K6free1}
For every $n$-vertex $K_6$-free graph $G$, $f_B(G) \le 4403n^2/15000$.
\end{THM}
\begin{Proof}
There are three possible cases:

\textbf{Case 1:} $b(G) \le 49n/100$.  By Lemma~\ref{l:PuleoTriIndependentSet1}
and Theorem~\ref{t:K6freeBipartite2}, we have
\begin{eqnarray*}
f_B(G) = \alpha_{1} (G) + \tau_{B} (G)
      \le \frac{nb(G)}{4} + \frac{17n^2}{100}
      \le \frac{49n^2}{400} + \frac{17n^2}{100}
      = \frac{117n^2}{400}
      < \frac{4403n^2}{15000}.
\end{eqnarray*}

\textbf{Case 2:} $49n/100 < b(G) \le 7n/10$.

By Lemma~\ref{l:PuleoTriIndependentSet1} and
Corollary~\ref{c:K6freeBipartite4} $(b)$, we have
\begin{eqnarray}
f_B(G) = \alpha_{1} (G) + \tau_{B} (G) \le \frac{nb(G)}{4} + \frac{17{(n-b(G))}^2}{150} + \frac{e}{3}.
\end{eqnarray}

By Lemma~\ref{l:PuleoTriIndependentSet2} and
Corollary~\ref{c:K6freeBipartite4} $(b)$, we have
\begin{eqnarray}
f_B(G) = \alpha_{1} (G) + \tau_{B} (G)
\nonumber &\le& \frac{n^2}{2} - e + \frac{17{(n-b(G))}^2}{150} + \frac{e}{3} \\
&=& \frac{n^2}{2} + \frac{17{(n-b(G))}^2}{150} - \frac{2e}{3}.
\end{eqnarray}

Multiplying inequality (1) by $2/3$ and inequality (2) by $1/3$, and adding them together, we have
\begin{eqnarray*}
f_B(G) \le \frac{nb(G)}{6} + \frac{n^2}{6} + \frac{17{(n-b(G))}^2}{150}
= \frac{1}{150} \left( 17b(G)^2 - 9nb(G) + 42n^2 \right).
\end{eqnarray*}
Since the function $g(x) = \left( 17x^2 - 9nx + 42n^2 \right)/150$ is increasing
in the interval $49n/100 < x \le 7n/10$, it follows that
\begin{eqnarray*}
f_B(G) \le g\hspace{-1mm}\left( \frac{7n}{10} \right)
= \frac{4403n^2}{15000}.
\end{eqnarray*}

\textbf{Case 3:} $b(G) > 7n/10$.
It is easy to verify that
$(-7b(G)^2+4nb(G)+3n^2)/20 \ge (-32b(G)^2+15nb(G)+17n^2)/100$
in this case.  So by Corollary~\ref{c:K6freeBipartite4} $(a)$,
we have $\tau_{B} (G) \le (-7b(G)^2+4nb(G)+3n^2)/20$.
Again, by Lemma~\ref{l:PuleoTriIndependentSet1} we have
$\alpha_{1} (G) \le nb(G)/4$.  Thus,
\begin{eqnarray*}
\noindent f_B(G) &=& \alpha_{1} (G) + \tau_{B} (G)\\
      &\le& \frac{nb(G)}{4} + \frac{-7b(G)^2+4nb(G)+3n^2}{20}\\
      &=& \frac{1}{20} \left( -7b(G)^2 + 9nb(G) + 3n^2 \right).
\end{eqnarray*}
The function $h(x) = (-7x^2 + 9nx + 3n^2)/20$ is
decreasing in the interval $x \ge 7n/10$.
So in this case we have
\begin{eqnarray*}
f_B(G) \le h\hspace{-1mm}\left( \frac{7n}{10} \right) = \frac{587n^2}{2000} < \frac{4403n^2}{15000}.
\end{eqnarray*}
\end{Proof}

\begin{THM} \label{t:Main2}
For every $n$-vertex graph $G$, $f_B(G) \le 4403n^2/15000$.
\end{THM}
\begin{Proof} It is easy to verify the conclusion for small $n$.
Now assume to the contrary that $G$ is a minimal counterexample.
By Theorem~\ref{t:K6free1}, $\omega (G) \ge 6$.
However, by Lemma~\ref{l:CliqueLemma} we have
$\omega (G) < 1/(4\times \frac{4403}{15000} - 1) = \frac{3750}{653} < 6$,
a contradiction.
\end{Proof}

\section{Acknowledgements}

The author would like to thank Wayne Goddard for
making many helpful comments and suggestions.
He also would like to thank Gregory Puleo for
pointing out an error in an earlier version of the paper.


\bibliographystyle{plain}
\bibliography{Bipartite}

\end{document}